\documentclass{amsart}

\usepackage[leqno]{amsmath}
\usepackage{latexsym,amsfonts,amssymb}
\usepackage[hidelinks]{hyperref}

\renewcommand{\eqref}[1]{(\ref{eq:#1})}
\def\urltilda{\kern -.15em\lower .7ex\hbox{\~{}}\kern .04em} 
\newcommand{\pff}[1]{\operatorname{\mathcal{P}}[#1]}
\newcommand{\sgn}[2]{\operatorname{sgn}%
  \big(\begin{smallmatrix}#1\\#2 \end{smallmatrix}\big)}
\newcommand{\Sign}[1]{(-1)^{#1}}
\newcommand{\dett}[2]{\mathcal{D}(#1;#2)}
 \newcommand{\Pf}[2][]{\operatorname{Pf}_{#1}(#2)}
 \newcommand{\set}[2][\mspace{1mu}]{\{#1 #2 #1\}}
\newcommand{\qtext}[1]{\quad\text{#1}\quad}
\newcommand{\qqtext}[1]{\qquad\text{#1}\qquad}
\newcommand{\qand}{\qtext{and}}
\newcommand{\qqand}{\qqtext{and}}
\newcommand{\deq}{\:=\:}
\newcommand{\gra}{\alpha}

\newcommand{\grr}{\rho}
\newcommand{\gro}{\omega}
\newcommand{\grs}{\sigma}

\theoremstyle{plain}
\newtheorem{theorem}[]{Theorem}
\numberwithin{equation}{subsection}

\begin{document}

\title[Minors of a skew symmetric matrix: A combinatorial
approach]{Minors of a skew symmetric matrix: \\A combinatorial
  approach}

\author[L.\,W. Christensen]{Lars Winther Christensen}

\address{Texas Tech University, Lubbock, TX 79409, U.S.A.}

\email{lars.w.christensen@ttu.edu}

\urladdr{http://www.math.ttu.edu/\urltilda lchriste}

\author[O. Veliche]{Oana Veliche}

\address{Northeastern University, Boston, MA~02115, U.S.A.}

\email{o.veliche@northeastern.edu}

\urladdr{https://web.northeastern.edu/oveliche}

\author[J. Weyman]{Jerzy Weyman}

\address{University of Connecticut, Storrs, CT~06269, U.S.A.}
\email{jerzy.weyman@uconn.edu}

\urladdr{http://www.math.uconn.edu/\urltilda weyman}

\thanks{L.W.C.\ was partly supported by Simons Foundation
  collaboration grant 428308, and J.W.\ was partly supported by NSF
  DMS grant 1802067.}

\date{29 July 2020}



\begin{abstract}
  We use Knuth's combinatorial approach to Pfaffians to reprove and
  clarify a century-old formula, due to Brill. It expresses arbitrary
  minors of a skew symmetric matrix in terms of Pfaffians.
\end{abstract}

\maketitle

\allowdisplaybreaks

\thispagestyle{empty}

\section{Brill's formula }

\noindent
In a paper \cite{DEK96} from 1996, Knuth took a combinatorial approach
to Pfaffians. It was immediately noticed that this approach
facilitates generalizations and simplified proofs of several known
identities involving Pfaffians; see for example Hamel \cite{AMH01}.

In this short note we apply the same approach to reprove a formula
that expresses arbitrary minors of a skew symmetric matrix in terms of
Pfaffians. The formula, which we here state as Theorem 1, first
appeared in a 1904 paper by Brill \cite{JBr04}. Our goal here is not
merely to give a short modern proof but also to clarify two aspects of
the formula: (1) To compute minors with Brill's formula
\cite[Eqn.~(1)]{JBr04} one needs to apply a sign convention explained
in a footnote; we integrate this sign into the formula. (2) As Brill
states right at the top of his paper, one can reduce to the case of
minors with disjoint row and column sets, and his formula deals only
with that case. We build that reduction into our proof and arrive at a
formula that also holds for minors with overlapping row and column
sets.

Bibliometric data suggests that Brill's formula may have gone largely
unnoticed. For us it turned out to be key to clarify a computation
in commutative algebra \cite{CVW-4}.

\subsection{Pfaffians following Knuth}
Let $T = (t_{ij})$ be an $n \times n$ skew symmetric matrix with
entries in some commutative ring. Assume that $T$ has zeros on the
diagonal; this is, of course, automatic if the characteristic of the
ring is not $2$. Set $\pff{ij} = t_{ij}$ for
$i,j \in \set{1,\ldots,n}$ and extend $\mathcal{P}$ to a function on
words in letters from the set $\set{1,\ldots,n}$ as follows:
\begin{equation*}
  \pff{i_1\ldots i_{m}} \deq
  \begin{cases}
    0 & \text{ if $m$ is odd} \\
    \sum \sgn{i_1\ldots i_{2k}}{j_1 \ldots j_{2k}}\pff{j_1j_2} \cdots
    \pff{j_{2k-1}j_{2k}} & \text{ if $m = 2k$ is even}
  \end{cases}
\end{equation*}
where the sum is over all partitions of $\set{i_1,\ldots,i_{2k}}$ into
$k$ subsets of cardinality $2$. The order of the two elements in each
subset is irrelevant as the difference in sign
$\pff{jj'} = - \pff{j'j}$ is offset by a change of sign of the
permutation; see \cite[Section~0]{DEK96}. The value of $\mathcal{P}$ on the
empty word is by convention $1$, and the value of $\mathcal{P}$ on a
word with a repeated letter is $0$. The latter is a convention in
characteristic $2$ and otherwise automatic.

For subsets $R$ and $S$ of $\set{1,\ldots,n}$ we write $T[R;S]$ for
the submatrix of $T$ obtained by taking the rows indexed by $R$ and
the columns indexed by $S$. The function $\mathcal{P}$ computes the
Pfaffians of skew symmetric submatrices of $T$. Indeed, for a subset
$R \subseteq \set{1,\ldots,n}$ with elements $r_1 < \cdots < r_m$ one
has
\begin{equation*}
  \Pf{T[R;R]} \deq \pff{r_1\ldots r_{m}} \:.
\end{equation*}
With this approach to Pfaffians, Knuth~\cite{DEK96} gives elegant
proofs for several classic formulas, generalizing them and clarifying
sign conventions in the process through the introduction of
permutation signs. For a word $\gra$ in letters from
$\set{1,\ldots,n}$ and $a \in \gra$ the following variation on the
classic Laplacian expansion is obtained as a special case of a formula
ascribed to Tanner \cite{HWT78},
\begin{equation}
  \label{eq:DEK}
  \pff{\gra} \deq \sum_{x\in\gra}
  \sgn{\gra}{ax(\gra\setminus ax)}
  \pff{ax}\pff{\gra\setminus ax} \:;
\end{equation}
see \cite[Eqn.~(2.0)]{DEK96}.  The introduction of the permutation
sign is what facilitates our statement and proof of Brill's formula.

\begin{theorem}
  Let $T$ be an $n \times n$ skew symmetric matrix. Let $R$ and $S$ be
  subsets of the set $\set{1,\ldots,n}$ with elements
  $r_1<\cdots <r_m$ and $s_1<\cdots <s_m$. With
  \begin{equation*}
    \grr \deq r_1\ldots r_m \qand
    \grs = s_1\ldots s_m
  \end{equation*}
  the next equality holds:
  \begin{equation*}
    \det(T[R;S]) 
    \deq \Sign{\left\lfloor \frac{m}{2} \right\rfloor}
    \sum_{k = 0}^{\left\lfloor
        \frac{m}{2} \right\rfloor}
    \mspace{6mu}\Sign{k}\mspace{-6mu}\sum_{\substack{|\omega| = 2k \\
        \gro \subseteq \grr}}
    \sgn{\grr}{\gro(\grr\setminus\gro)}\pff{\omega}
    \pff{(\grr\setminus\omega)\grs}\:.
  \end{equation*}
\end{theorem}

\noindent Notice that only subwords $\omega$ of $\rho$ that contain
$\grr \cap \grs$ contribute to the sum, otherwise the word
$(\grr\setminus\omega)\grs$ has a repeated letter. The quantity
$\pff{\gro}$ is the Pfaffian of the submatrix
$\Pf{T[\set{\gro};\set{\gro}]}$, where $\set{\gro}$ is the subset of
$\set{1,\ldots,n}$ whose elements are the letters in
$\gro$. Similarly, for a word $\gro$ that contains $\grr \cap \grs$,
the quantity $\pff{(\grr\setminus\omega)\grs}$ equals
$\sgn{\gra}{(\gra\setminus\grs)\grs} \Pf{T[\set{\gra};\set{\gra}]}$
where the word $\gra$ has the letters from $(\grr\setminus\omega)\grs$
written in increasing order.

\section{Our proof}

\noindent
For subsets $R$ and $S$ as in the statement set
\begin{equation*}
  \dett{\grr}{\grs} = \det(T[R;S]) \:.
\end{equation*}
For a word $\gra$ in letters from $\set{1,\ldots,n}$ we write
$\overline{\gra}$ for the word that has the letters from $\gra$
written in increasing order. In the words $\grr$ and $\grs$ and their
subwords the letters already appear in increasing order, but the bar
notation comes in handy for concatenated words such as $\grr\grs$.

\subsection{Disjoint row and column sets}
We first assume that $\grr$ and $\grs$ are disjoint subwords of
$1\ldots n$. Under this assumption
$\pff{(\grr\setminus\omega)\grs} = \pff{\grr\grs\setminus\omega}$
holds for every subword $\gro \subseteq \grr$, and in the balance of
the section we use the simpler notation
$\pff{\grr\grs\setminus\omega}$.  We proceed by induction on $m$.  For
$m=1$ the formula holds as one has
$\det(T[\set{r_1};\set{s_1}]) = t_{r_1s_1} = \pff{r_1s_1}$. For
$m = 2$ expansion of the determinant along row $r_1$ yields
\begin{equation*}
  \dett{\grr}{\grs} = \pff{r_1s_1}\pff{r_2s_2} - \pff{r_1s_2}\pff{r_2s_1} \:.
\end{equation*}
With $\gra = \grr\grs$ and $a = r_1$ the formula \eqref{DEK}
reads
\begin{align*}
  \pff{\grr\grs} 
  & \deq \sum_{x\in\grr\grs}
    \sgn{\grr\grs}{r_1x(\grr\grs\setminus r_1x)}
    \pff{r_1x}\pff{\grr\grs\setminus r_1x}
  \\
  & \deq \pff{r_1r_2}\pff{s_1s_2} -
    \pff{r_1s_1}\pff{r_2s_2} + \pff{r_1s_2}\pff{r_2s_1}
  \\
  & \deq \pff{\grr}\pff{\grs} - \dett{\grr}{\grs} \:,
\end{align*}
which can be rewritten in the desired form:
\begin{equation*}
  \dett{\grr}{\grs} \deq -( \pff{\grr\grs} - \pff{\grr}\pff{\grs} ) \:.
\end{equation*}

Before we move on to the induction step we record another consequence
of \eqref{DEK}:
\begin{equation}
  \label{eq:rso}
  \begin{aligned}
    \pff{\grr\grs\setminus\gro} 
    & \deq \sum_{r\in\grr\setminus r_1\gro}
    \sgn{\grr\grs\setminus\gro}{r_1r(\grr\grs\setminus\gro r_1r)}
    \pff{r_1r}\pff{\grr\grs\setminus \gro r_1r} \\
    & \hspace{4pc} {} + \sum_{s\in\grs}
    \sgn{\grr\grs\setminus\gro}{r_1s(\grr\grs\setminus\gro r_1s)}
    \pff{r_1s}\pff{\grr\grs\setminus \gro r_1s}
    \\
    & \deq \sum_{r\in\grr\setminus r_1\gro}
    \sgn{\grr\setminus\gro}{r_1r(\grr\setminus\gro r_1r)}
    \pff{r_1r}\pff{\grr\grs\setminus \gro r_1r} \\
    & \hspace{4pc} {} + \Sign{|\grr \setminus \gro| - 1}
    \sum_{s\in\grs} \sgn{\grs}{s(\grs\setminus s)}
    \pff{r_1s}\pff{\grr\grs\setminus \gro r_1s} \:.
  \end{aligned}
\end{equation}

Now let $m \ge 2$ and $|\grr| = m+1 = |\grs|$. In the next computation
the first equality is the expansion of the determinant
$\dett{\grr}{\grs}$ along row $r_1$. The second equality follows from
the induction hypothesis, and the third is obtained by changing the
order of summation. The fourth equality follows from \eqref{rso}. The
fifth equality uses that $|\gro|$ is even and that
$|\grr\setminus \gro| -1 = m - 2k \equiv m \mod 2$ holds.
\begin{align*}
  \dett{\grr}{\grs} 
  & \deq \sum_{s\in\grs}\sgn{\grs}{s\grs\setminus s}\pff{r_1s}
    \,\dett{\grr\setminus r_1}{\grs\setminus s}
  \\
  & \deq \sum_{s\in\grs}\sgn{\grs}{s\grs\setminus s}\pff{r_1s}  \\
  & \hspace{4pc} \cdot \Biggl(
    \Sign{\left\lfloor \frac{m}{2} \right\rfloor}
    \sum_{k = 0}^{\left\lfloor \frac{m}{2} \right\rfloor}
    \mspace{6mu}\Sign{k}\mspace{-6mu} \sum_{\substack{
    |\omega| = 2k \\ \scriptscriptstyle \gro \subseteq \grr\setminus r_1}}
  \mspace{-4mu}\sgn{\grr\setminus r_1}{\gro(\grr\setminus r_1\gro)}
  \pff{\omega}\pff{\grr\grs\setminus r_1\omega s}
  \Biggr)
  \\
  & \deq \Sign{\left\lfloor \frac{m}{2} \right\rfloor}
    \sum_{k =0}^{\left\lfloor \frac{m}{2} \right\rfloor}
    \mspace{6mu}\Sign{k}\mspace{-6mu} \sum_{\substack{|\omega| = 2k \\
  \scriptscriptstyle \gro \subseteq \grr\setminus r_1}}
  \mspace{-4mu}\sgn{\grr\setminus r_1}{\gro(\grr\setminus r_1\gro)} \pff{\omega} 
  \\ & \hspace{13pc} \cdot \Bigr(
       \sum_{s\in\grs}\sgn{\grs}{s\grs\setminus s}
       \pff{r_1s} \pff{\grr\grs\setminus r_1\omega s}
       \Bigl) 
  \\
  & \deq \Sign{\left\lfloor \frac{m}{2} \right\rfloor}
    \sum_{k =0}^{\left\lfloor \frac{m}{2} \right\rfloor}
    \mspace{6mu}\Sign{k}\mspace{-6mu} \sum_{\substack{|\omega| = 2k \\
  \scriptscriptstyle \gro \subseteq \grr\setminus r_1}}
  \mspace{-4mu}\sgn{\grr\setminus r_1}{\gro(\grr\setminus r_1\gro)} \pff{\omega} 
  \\ & \hspace{2pc}
       \cdot \Sign{|\grr\setminus\gro|-1}
       \Big(\pff{\grr\grs\setminus\gro}
       \mspace{4mu} - \mspace{-8mu}\sum_{r\in\grr\setminus\gro r_1}
       \mspace{-9mu}\sgn{\grr\setminus\gro}{r_1r(\grr\setminus\gro r_1r)}
       \pff{r_1r}\pff{\grr\grs\setminus\gro r_1r} \Big)
  \\
  & \deq \Sign{\left\lfloor \frac{m}{2} \right\rfloor + m}
    \Biggl(
    \sum_{k=0}^{\left\lfloor \frac{m}{2} \right\rfloor} 
    \mspace{6mu}\Sign{k}\mspace{-6mu} \sum_{\substack{
    |\omega| = 2k \\ \scriptscriptstyle \gro \subseteq \grr\setminus r_1}}
  \sgn{\grr}{\gro(\grr\setminus\gro)}
  \pff{\omega} \pff{\grr\grs\setminus\gro} - 
  \sum_{k=0}^{\left\lfloor
  \frac{m}{2} \right\rfloor}
  \Sign{k}  \\
  & \hspace{2pc} \cdot \mspace{-4mu} \sum_{\substack{|\omega| = 2k \\
  \scriptscriptstyle \gro \subseteq \grr\setminus r_1}}
  \mspace{-4mu}\sgn{\grr}{\gro(\grr\setminus \gro)}
  \pff{\omega} \mspace{-4mu}
  \sum_{r\in\grr\setminus\gro r_1}
  \mspace{-9mu}\sgn{\grr\setminus\gro}{r_1r(\grr\setminus\gro r_1r)}
  \pff{r_1r}\pff{\grr\grs\setminus\gro r_1r}
  \Bigg)
  \\ 
  & \deq \Sign{\left\lfloor \frac{m+1}{2} \right\rfloor}
    \Bigg(
    \sum_{k=0}^{\left\lfloor \frac{m}{2} \right\rfloor} 
    \mspace{6mu}\Sign{k}\mspace{-6mu} \sum_{\substack{
    |\omega| = 2k \\ \scriptscriptstyle \gro \subseteq \grr\setminus r_1}}
  \mspace{-4mu}\sgn{\grr}{\gro(\grr\setminus\gro)}
  \pff{\omega} \pff{\grr\grs\setminus\gro} + 
  \sum_{k=0}^{\left\lfloor
  \frac{m}{2} \right\rfloor}
  \Sign{k+1} \\ &\hspace{3pc} \cdot \mspace{-3mu}  \sum_{\substack{
                  |\omega| = 2k \\ \scriptscriptstyle \gro \subseteq \grr\setminus r_1}}
  \sum_{r\in\grr\setminus\gro r_1}
  \mspace{-9mu}\sgn{\grr}{\gro r_1r(\grr\setminus \gro r_1r)}
  \pff{\omega} \pff{r_1r}\pff{\grr\grs\setminus\gro r_1r}
  \Bigg) \:.
\end{align*}
The next step is to simplify the last line in the display above. The
first equality in the following computation holds as $|\gro|$ is
even. The third equality follows by substituting $\gro'$ for the word
$\overline{r_1r\gro}$ and noticing that one then has
\begin{equation*}
  \sgn{\grr}{r_1r\gro (\grr\setminus r_1r\gro)} \deq
  \sgn{\grr}{\gro'(\grr\setminus \gro')}
  \sgn{\gro'}{r_1r\gro }
  \deq
  \sgn{\grr}{\gro'(\grr\setminus \gro')}
  \sgn{\gro'}{r_1r(\gro'\setminus r_1r)} \:.
\end{equation*}
The last equality follows from \eqref{DEK} applied with
$\gra = \gro'$ and $a = r_1$.
\begin{align*}
  {} & \sum_{\substack{|\omega| = 2k \\ \scriptscriptstyle \gro \subseteq \grr\setminus r_1}}
  \mspace{3mu}\sum_{r\in\grr\setminus\gro r_1}
  \mspace{-4mu}\sgn{\grr}{\gro r_1r(\grr\setminus \gro r_1r)}
  \pff{\omega} \pff{r_1r}\pff{\grr\grs\setminus\gro r_1r}
  \\
     &\deq \sum_{\substack{
       |\omega| = 2k \\ \scriptscriptstyle \gro \subseteq \grr\setminus r_1}}
  \mspace{3mu}\sum_{r\in\grr\setminus\gro r_1}
  \sgn{\grr}{r_1r \gro(\grr\setminus r_1r\gro)}
  \pff{\omega} \pff{r_1r}\pff{\grr\grs\setminus r_1r\gro}
  \\
     & \deq \sum_{\substack{|\omega'| = 2k+2 \\ \scriptscriptstyle r_1 \in \gro' \subseteq \grr}}
  \mspace{3mu} \sum_{r\in\gro'\setminus r_1}
  \mspace{-4mu}  \sgn{\grr}{\gro'(\grr\setminus \gro')}
  \sgn{\gro'}{r_1r(\gro'\setminus r_1r)}
  \pff{\omega'\setminus r_1r} \pff{r_1r}\pff{\grr\grs\setminus\gro'}
  \\
     & \deq \sum_{\substack{
       |\omega'| = 2k+2 \\ \scriptscriptstyle r_1 \in \gro' \subseteq \grr}}
  \mspace{-4mu}\sgn{\grr}{\gro'(\grr\setminus \gro')}
  \pff{\grr\grs\setminus\gro'} 
  \Big( 
  \sum_{r\in\gro'\setminus r_1}
  \mspace{-4mu}\sgn{\gro'}{r_1r(\gro'\setminus r_1r)}
  \pff{r_1r}\pff{\omega'\setminus r_1r} 
  \Big)
  \\
     & \deq \sum_{\substack{ |\omega'| = 2(k+1) \\
  \scriptscriptstyle r_1 \in \gro' \subseteq \grr}}
  \mspace{-4mu} \sgn{\grr}{\gro'(\grr\grs\setminus \gro')}
  \pff{\grr\setminus\gro'} \pff{\gro'} \:.
\end{align*}
Substituting this into the computation above one gets the desired
equality.
\begin{align*}
  \dett{\grr}{\grs} 
  & \deq \Sign{\left\lfloor \frac{m+1}{2} \right\rfloor}
    \Bigg(
    \sum_{k=0}^{\left\lfloor \frac{m}{2} \right\rfloor} 
    \mspace{6mu}\Sign{k}\mspace{-6mu} \sum_{\substack{
    |\omega| = 2k \\ \scriptscriptstyle \gro \subseteq \grr\setminus r_1}}
  \mspace{-3mu}\sgn{\grr}{\gro(\grr\setminus\gro)}
  \pff{\omega} \pff{\grr\grs\setminus\gro}  \\
  & \hspace{4pc} {} + \sum_{k=0}^{\left\lfloor
    \frac{m}{2} \right\rfloor}
    \mspace{6mu}\Sign{k+1}\mspace{-12mu}
    \sum_{\substack{ |\omega'| = 2(k+1) \\
  \scriptscriptstyle r_1 \in \gro' \subseteq \grr}}
  \mspace{-4mu}\sgn{\grr}{\gro'(\grr\setminus \gro')}
  \pff{\grr\grs\setminus\gro'} \pff{\gro'} \Bigg)
  \\
  & \deq \Sign{\left\lfloor \frac{m+1}{2} \right\rfloor}
    \sum_{k=0}^{\left\lfloor \frac{m+1}{2} \right\rfloor} 
    \Sign{k} \sum_{\substack{ |\omega| = 2k \\
  \scriptscriptstyle \gro \subseteq \grr}}
  \sgn{\grr}{\gro(\grr\setminus\gro)}
  \pff{\omega} \pff{\grr\grs\setminus\gro} \:.
\end{align*}
This proves the asserted formula in the special case where $R$ and $S$
are disjoint, and it remains to reduce the general case to this one.

\subsection{Overlapping row and column sets}
First notice that the submatrix $T[R;S]$ agrees with a submatrix of
the $2n \times 2n$ skew symmetric matrix $T'$ obtained from $T$ by
repeating all entries, horizontally and vertically. For example, with
{\small \begin{equation*} T \deq
    \begin{pmatrix}
      0 & t_{12} & t_{13} \\
      - t_{12} & 0 &  t_{23} \\
      - t_{13} & - t_{23} & 0
    \end{pmatrix}
    \qand
    T' \deq
    \begin{pmatrix}
      0 & 0 & t_{12} & t_{12} & t_{13} & t_{13} \\
      0 & 0 & t_{12} & t_{12} & t_{13} & t_{13} \\
      - t_{12} & - t_{12} & 0 & 0 &  t_{23} &  t_{23} \\
      - t_{12} & - t_{12} & 0 & 0 &  t_{23} &  t_{23} \\
      - t_{13} & - t_{13} & - t_{23} & - t_{23} & 0 & 0 \\
      - t_{13} & - t_{13} & - t_{23} & - t_{23} & 0 & 0
    \end{pmatrix}
  \end{equation*} }
one has $T[\set{1,2};\set{2,3}] = T'[\set{1,3};\set{4,6}]$. In
general, one has 
\begin{gather*}
  T[R;S] \deq T'[R';S'] \intertext{for} R' \deq \set{2r_1-1,\ldots,
    2r_m-1} \qqand S' = \set{2s_1,\ldots, 2s_m} \:.
\end{gather*}
One now has
\begin{align*}
  \det(T[R;S]) & \deq \det(T'[R';S']) \\
               & \deq \Sign{\left\lfloor \frac{m}{2} \right\rfloor}
                 \sum_{k=0}^{\left\lfloor
                 \frac{m}{2} \right\rfloor}
                 \mspace{6mu}\Sign{k}\mspace{-6mu}\sum_{\substack{ |\omega'| = 2k \\
  \scriptscriptstyle \gro' \subseteq \grr'}}
  \sgn{\grr'}{\gro'(\grr'\setminus\gro')}\pff{\gro'}
  \pff{\grr'\grs'\setminus\gro'}
\end{align*}
where the second equality holds by the case already
established. Indeed, the numbers in $R'$ are odd and the numbers in
$S'$ are even.  Notice that for every $k$ there is a one-to-one
correspondence, given by the relation $r_i' = 2r_i-1$, between
subwords $\gro'\subseteq \grr'$ of length $2k$ and subwords $\gro$ of
$\grr$ of length $2k$, and one has $\pff{\gro'} = \pff{\gro}$ for
corresponding subwords. Let $U$ denote the submatrix of $T'[R',S']$
obtained by removing the rows and columns whose indices, all odd,
appear in the word $\gro'$. One has
\begin{equation*}
  \pff{\grr'\grs'\setminus\omega'} = \sgn{\overline{\grr'\grs'\setminus\omega'}}{\grr'\grs'\setminus\omega'} \Pf{U} \:.
\end{equation*}
If an element $x \in R \cap S$ does not appear in $\gro$, then $2x-1$
is not in $\gro'$, so $U$ contains the identical rows/columns $2x-1$
and $2x$ from $T'$, whence
$\pff{\grr'\grs'\setminus\omega'}=0=\pff{(\grr\setminus\omega)\grs}$
holds. On the other hand, if every element of $R \cap S$ appears in
$\gro$, then the matrix $U$ is a submatrix of
$T[R\setminus R \cap S; S\setminus R \cap S]$ with Pfaffian
$\sgn{\overline{(\grr\setminus\omega)\grs}}{(\grr\setminus\omega)\grs}\pff{(\grr\setminus\omega)\grs}$. As
one has $\grr'\grs'\setminus\omega' = (\grr'\setminus\omega')\grs'$
the permutations
$\big(\begin{smallmatrix}\overline{(\grr\setminus\omega)\grs}\\
  (\grr\setminus\omega)\grs \end{smallmatrix}\big)$ and
$\big(\begin{smallmatrix}\overline{\grr'\grs'\setminus\omega'}\\
  \grr'\grs'\setminus\omega' \end{smallmatrix}\big)$ have the same
sign, and the proof is complete.

  \providecommand{\MR}[1]{\mbox{\href{http://www.ams.org/mathscinet-getitem?mr=#1}{#1}}}
  \renewcommand{\MR}[1]{\mbox{\href{http://www.ams.org/mathscinet-getitem?mr=#1}{#1}}}
  \providecommand{\arxiv}[2][AC]{\mbox{\href{http://arxiv.org/abs/#2}{\sf
  arXiv:#2 [math.#1]}}} \def\cprime{$'$}
\providecommand{\bysame}{\leavevmode\hbox to3em{\hrulefill}\thinspace}
\providecommand{\MR}{\relax\ifhmode\unskip\space\fi MR }
\providecommand{\MRhref}[2]{%
  \href{http://www.ams.org/mathscinet-getitem?mr=#1}{#2}
}
\providecommand{\href}[2]{#2}

\end{document}